# Consistency and robustness of kernel-based regression in convex risk minimization

ANDREAS CHRISTMANN[1] and INGO STEINWART[2]

[1]*Department of Mathematics, Vrije Universiteit Brussel, Brussels, Belgium.*
*E-mail: Andreas.Christmann@vub.ac.be*
[2]*Los Alamos National Laboratory, CCS-3, Los Alamos, New Mexico, USA.*
*E-mail: ingo@lanl.gov*

We investigate statistical properties for a broad class of modern kernel-based regression (KBR) methods. These kernel methods were developed during the last decade and are inspired by convex risk minimization in infinite-dimensional Hilbert spaces. One leading example is support vector regression. We first describe the relationship between the loss function $L$ of the KBR method and the tail of the response variable. We then establish the $L$-risk consistency for KBR which gives the mathematical justification for the statement that these methods are able to "learn". Then we consider robustness properties of such kernel methods. In particular, our results allow us to choose the loss function and the kernel to obtain computationally tractable and consistent KBR methods that have bounded influence functions. Furthermore, bounds for the bias and for the sensitivity curve, which is a finite sample version of the influence function, are developed, and the relationship between KBR and classical $M$ estimators is discussed.

*Keywords:* consistency; convex risk minimization; influence function; nonparametric regression; robustness; sensitivity curve; support vector regression

## 1. Introduction

The goal in nonparametric regression is to estimate a functional relationship between an $\mathbb{R}^d$-valued input random variable $X$ and an $\mathbb{R}$-valued output random variable $Y$, under the assumption that the joint distribution P of $(X, Y)$ is (almost) completely *unknown*. To solve this problem, one typically assumes a set of observations $(x_i, y_i)$ from independent and identically distributed (i.i.d.) random variables $(X_i, Y_i)$, $i = 1, \ldots, n$, which all have the distribution P. Informally, the aim is to build a predictor $f: \mathbb{R}^d \to \mathbb{R}$ on the basis of these observations such that $f(X)$ is a "good" approximation of $Y$. To formalize this aim, one assumes a *loss* $L$, that is, a continuous function $L: \mathbb{R} \times \mathbb{R} \to [0, \infty)$, that assesses the quality of a prediction $f(x)$ for an observed output $y$ by $L(y, f(x))$. Here it is commonly assumed that the smaller $L(y, f(x))$ is, the better the prediction is. The quality of a predictor $f$ is then measured by the $L$ risk

$$\mathcal{R}_{L,\mathrm{P}}(f) := \mathbb{E}_{\mathrm{P}} L(Y, f(X)), \qquad (1)$$

that is, by the average loss obtained by predicting with $f$. Following the interpretation that a small loss is desired, one tries to find a predictor with risk close to the optimal risk $\mathcal{R}^*_{L,\mathrm{P}} := \inf\{\mathcal{R}_{L,\mathrm{P}}(f) | f : \mathbb{R}^d \to \mathbb{R} \text{ measurable}\}$. Traditionally, most research in nonparametric regression considered the least squares loss $L(y,t) := (y-t)^2$, mainly because it *"simplifies the mathematical treatment"* and *"leads naturally to estimates which can be computed rapidly"* as Györfi *et al.* [12], page 2, wrote. However, from a practical point of view, there are situations in which a different loss is more appropriate, for example:

**Regression problem not described by least squares loss.** It is well known that the least squares risk is minimized by the conditional mean $\mathbb{E}_{\mathrm{P}}(Y|X=x)$. However, in many situations one is actually not interested in this mean, but in, for example, the conditional median instead. Now recall that the conditional median is the minimizer of $\mathcal{R}_{L,\mathrm{P}}$, where $L$ is the absolute value loss, that is, $L(y,t) := |y-t|$, and the same statement holds for conditional quantiles if one replaces the absolute value loss by an *asymmetric* variant known as the pinball loss; see Steinwart [24] for approximate minimizers and Christmann and Steinwart [6] for kernel-based quantile regression.

**Surrogate losses.** If the conditional distributions of $Y|X=x$ are known to be symmetric, then basically all loss functions of the form $L(y,t) = l(y-t)$, where $l : \mathbb{R} \to [0,\infty)$ is convex, symmetric and has its only minimum at 0, can be used to estimate the conditional mean; see Steinwart [24]. In this case, a less steep surrogate such as absolute value loss, Huber's loss (given by $l(r) = r^2$ if $|r| \leq c$, and $l(r) = c|r| - c^2/2$ if $|r| > c$ for some $c \in (0,\infty)$), or logistic loss (given by $l(r) = -\log\{4\Lambda(r)[1 - \Lambda(r)]\}$, where $\Lambda(r) := 1/(1+e^{-r})$) may be more suitable if one expects outliers in the $y$ direction.

**Algorithmic aspects.** If the goal is to estimate the conditional median, then Vapnik's $\epsilon$-insensitive loss, given by $l(r) = \max\{|r| - \epsilon, 0\}$, $\epsilon \in (0,\infty)$, promises algorithmic advantages in terms of sparseness compared to the absolute loss when it is used in the kernel-based regression method described below. This is especially important for large data sets, which commonly occur in data mining.

One way to build a nonparametric predictor $f$ is to use kernel-based regression (KBR), which finds a minimizer $\hat{f}_{n,\lambda}$ of the regularized empirical risk

$$\frac{1}{n} \sum_{i=1}^{n} L(y_i, f(x_i)) + \lambda \|f\|_H^2, \qquad (2)$$

where $\lambda > 0$ is a regularization parameter to reduce the danger of overfitting, $H$ is a reproducing kernel Hilbert space (RKHS) of a kernel $k : X \times X \to \mathbb{R}$ and $L$ is a *convex* loss in the sense that $L(y, \cdot) : \mathbb{R} \to [0, \infty)$ is convex for all $y \in Y$. Because (2) is strictly convex in $f$, the minimizer $\hat{f}_{n,\lambda}$ is uniquely determined and a simple gradient descent algorithm can be used to find $\hat{f}_{n,\lambda}$. However, for specific losses such as the least squares or the $\epsilon$-insensitive loss, other, more efficient algorithmic approaches are used in practice; see Wahba [29], Vapnik [28], Schölkopf and Smola [20] and Suykens *et al.* [26].

Of course, when using KBR, one natural question is whether the risk $\mathcal{R}_{L,\mathrm{P}}(\hat{f}_{n,\lambda})$ actually tends to the minimal risk $\mathcal{R}^*_{L,\mathrm{P}}$ if $n \to \infty$. For example, if $Y$ is bounded, $H$ is rich in the sense of Definition 11 below and $\lambda = \lambda_n \to 0$ "sufficiently slowly", this question can be positively answered by current standard techniques based on uniform deviation inequalities. However, for *unbounded* $Y$, no such results presently exist. The first goal of this work is to close this gap (using Theorem 12).

Our second aim is to investigate robustness properties of KBR. To this end we consider

$$R^{\mathrm{reg}}_{L,\mathrm{P},\lambda}(f) := \mathbb{E}_{\mathrm{P}} L(Y, f(X)) + \lambda \|f\|_H^2, \qquad f \in H, \tag{3}$$

which can be interpreted as the infinite sample version of (2). In Section 4, we describe the influence of both the kernel and the loss function on the robustness of the (uniquely determined) minimizer $f_{\mathrm{P},\lambda}$ of (3). In particular, we establish the existence of the influence function for a broad class of KBR methods and present conditions under which the influence function is bounded. Here it turns out that, depending on the kernel and the loss function, some KBR methods are robust while others are not. Consequently, our results can help us to choose both quantities to obtain consistent KBR estimators with good robustness properties in a situation that allows such freedom. Moreover, it is interesting, but in some sense not surprising, that the robust KBR methods are exactly the ones that require the mildest tail conditions on $Y$ for consistency.

## 2. Some basics on losses, risks and kernel-based regression

In this section we first introduce some important concepts for loss functions that are used throughout this work. Thereafter, we investigate properties of their associated risks and discuss the interplay between growth behaviour of the loss functions and the tail of the response variable. Finally, we establish existence and stability results for the infinite-sample KBR methods given by (3). These results are needed to obtain both the consistency results in Section 3 and some of the robustness results in Section 4.

**Definition 1.** *Let $Y \subset \mathbb{R}$ be a non-empty closed subset and let $L: Y \times \mathbb{R} \to [0, \infty)$ be a measurable loss function. Then $L$ is called invariant if there exists a function $l: \mathbb{R} \to [0, \infty)$ with $l(0) = 0$ and $L(y, t) = l(y - t)$ for all $y \in Y$, $t \in \mathbb{R}$. Moreover, $L$ is called Lipschitz continuous if there exists a constant $c > 0$ such that*

$$|L(y, t) - L(y, t')| \leq c \cdot |t - t'| \tag{4}$$

*for all $y \in Y$, $t, t' \in \mathbb{R}$. In this case, we denote the smallest possible $c$ in (4) by $|L|_1$.*

Note that an invariant loss function $L$ is convex if and only if the corresponding $l$ is convex. Analogously, $L$ is Lipschitz continuous if and only if $l$ is Lipschitz continuous:

in this case, we have $|L|_1 = |l|_1$, where $|l|_1$ denotes the Lipschitz constant of $l$. Many popular loss functions for regression problems are convex, Lipschitz continuous and invariant. Three important examples are the Vapnik $\epsilon$-insensitive loss, the Huber loss and the logistic loss functions introduced in the Introduction. Moreover, the least squares loss function $L(y,t) = (y-t)^2$ is convex and invariant, but not Lipschitz continuous. The logistic loss function is a compromise between the other three loss functions: it is twice continuously differentiable with $L'' > 0$, which is true for the least squares loss function, and it increases approximately linearly if $|y-t|$ tends to infinity, which is true for Vapnik's and Huber's loss functions. These four loss functions are even symmetric, because $L(y,t) = L(t,y)$ for $y,t \in \mathbb{R}$. Asymmetric loss functions may be interesting in some applications where extremely skewed distributions occur, for example, analysis of claim sizes in insurance data (see Christmann [4]).

The growth behaviour of $L$ plays an important role in both consistency and robustness results. Hence we now introduce some basic concepts that describe the growth behaviour.

**Definition 2.** *Let $L : Y \times \mathbb{R} \to [0,\infty)$ be a loss function, let $a : Y \to [0,\infty)$ be a measurable function and let $p \in [0,\infty)$. We say that $L$ is a loss function of type $(a,p)$ if there exists a constant $c > 0$ such that $L(y,t) \leq c(a(y) + |t|^p + 1)$ for all $y \in Y$ and all $t \in \mathbb{R}$. We say that $L$ is of strong type $(a,p)$ if the first two partial derivatives $L' := \partial_2 L$ and $L := \partial_{22} L$ of $L$ with respect to the second argument of $L$ exist, and $L$, $L'$ and $L''$ are of $(a,p)$ type.*

For invariant loss functions, it turns out that there is an easy way to determine their type. To describe the corresponding results, we need the following definition.

**Definition 3.** *Let $L$ be an invariant loss function with corresponding function $l : \mathbb{R} \to \mathbb{R}$ and let $p \geq 0$. We say that $L$ is of upper order $p$ if there exists a constant $c > 0$ such that $l(r) \leq c(|r|^p + 1)$ for all $r \in \mathbb{R}$. Analogously, we say that $L$ is of lower order $p$ if there exists a constant $c > 0$ such that $l(r) \geq c(|r|^p - 1)$ for all $r \in \mathbb{R}$.*

Recalling that convex functions are locally Lipschitz continuous, we see that for invariant losses $L$, the corresponding $l$ is Lipschitz continuous on every interval $[-r,r]$. Consequently, $V(r) := |l_{|[-r,r]}|_1$, where $|l_{|[-r,r]}|_1$ denotes the Lipschitz constant of the restriction $l_{|[-r,r]}$ of $l$ onto $[-r,r]$ for $r \geq 0$, defines a non-decreasing function $V : [0,\infty) \to [0,\infty)$. We denote its symmetric extension also by $V$, so that we have $V(-r) = V(r)$ for all $r \in \mathbb{R}$. The next result lists some properties of invariant losses.

**Lemma 4.** *Let $L$ be an invariant loss function with corresponding $l : \mathbb{R} \to \mathbb{R}$ and $p \geq 0$.*

  (i) *If $L$ is convex and satisfies $\lim_{|r| \to \infty} l(r) = \infty$, then it is of lower order 1.*
  (ii) *If $L$ is Lipschitz continuous, then it is of upper order 1.*
  (iii) *If $L$ is convex, then for all $r > 0$ we have $V(r) \leq \frac{2}{r} \|l_{|[-2r,2r]}\|_\infty \leq 4V(2r)$.*
  (iv) *If $L$ is of upper order $p$, then $L$ is of type $(a,p)$ with $a(y) := |y|^p$, $y \in Y$.*

The proof of this lemma as well as the proofs of all following results can be found in the Appendix. With the help of Lemma 4, it is easy to see that the least squares loss function is of strong type $(y^2, 2)$. Furthermore, the logistic loss function is of strong type $(|y|, 1)$ because it is twice continuously differentiable with respect to its second variable and both derivatives are bounded, namely $|\partial_2 L(y,t)| \leq 1$ and $|\partial_{22} L(y,t)| \leq \frac{1}{2}$, $t \in \mathbb{R}$. The Huber and Vapnik loss functions are of upper and lower order 1 because they are convex and Lipschitz continuous; however, they are not of any strong type because they are not twice continuously differentiable.

Our next goal is to find a condition that ensures $\mathcal{R}_{L,P}(f) < \infty$. To this end we need the following definition, which for later purposes is formulated in a rather general way (see Brown and Pearcy [2] for signed measures).

**Definition 5.** *Let $\mu$ be a signed measure on $X \times Y$ with total variation $|\mu|$ and let $a: Y \to [0, \infty)$ be a measurable function. Then we write $|\mu|_a := \int_{X \times Y} a(y) \, d|\mu|(x,y)$. If $a(y) = |y|^p$ for some $p > 0$ and all $y \in Y$, we write $|\mu|_p := |\mu|_a$ whenever no confusion can arise. Finally, we write $|\mu|_0 := \|\mu\|_{\mathcal{M}}$, where $\|\mu\|_{\mathcal{M}}$ denotes the norm of total variation.*

We can now formulate the following two results investigating finite risks.

**Proposition 6.** *Let $L$ be an $(a,p)$-type loss function, let $P$ be a distribution on $X \times Y$ with $|P|_a < \infty$ and let $f: X \to \mathbb{R}$ be a function with $f \in L_p(P)$. Then we have $\mathcal{R}_{L,P}(f) < \infty$.*

**Lemma 7.** *Let $L$ be an invariant loss of lower order $p$, let $f: X \to \mathbb{R}$ be measurable and let $P$ be a distribution $X \times Y$ with $\mathcal{R}_{L,P}(f) < \infty$. Then $|P|_p < \infty$ if and only if $f \in L_p(P)$.*

If $L$ is an invariant loss function of lower and upper order $p$ and $P$ is a distribution with $|P|_p = \infty$, Lemma 7 shows $\mathcal{R}_{L,P}(f) = \infty$ for all $f \in L_p(P)$. This suggests that we may even have $\mathcal{R}_{L,P}(f) = \infty$ for *all* measurable $f: X \to Y$. However, this is in general not the case. For example, let $P_X$ be a distribution on $X$ and let $g: X \to \mathbb{R}$ be a measurable function with $g \notin L_p(P_X)$. Furthermore, let $P$ be the distribution on $X \times \mathbb{R}$ whose marginal distribution on $X$ is $P_X$ and whose conditional probability satisfies $P(Y = g(x)|x) = 1$. Then we have $|P|_p = \infty$, but $\mathcal{R}_{L,P}(g) = 0$.

Our next goal is to establish some preliminary results on (3). To this end, recall that the *canonical feature map* of a kernel $k$ with RKHS $H$ is defined by $\Phi(x) := k(\cdot, x)$, $x \in X$. Moreover, the *reproducing property* gives $f(x) = \langle f, k(\cdot, x) \rangle$ for all $f \in H$ and $x \in X$. Of special importance in terms of applications is the Gaussian radial basis function (RBF) kernel $k(x, x') = \exp(-\gamma \|x - x'\|^2)$, $\gamma > 0$, which is a universal kernel on every compact subset of $\mathbb{R}^d$; see Definition 11. This kernel is also *bounded* because $\|k\|_\infty = \sup\{\sqrt{k(x,x)} : x \in \mathbb{R}^d\} = 1$. Polynomial kernels $k(x, x') = (c + \langle x, x' \rangle)^m$, $m \geq 1$, $c \geq 0$, $x, x' \in \mathbb{R}$, are also popular in practice, but obviously they are neither universal nor bounded. Finally, recall that for bounded kernels, the canonical feature map satisfies $\|\Phi(x)\|_H \leq \|k\|_\infty$ for all $x \in \mathbb{R}^d$.

Let us now recall a result from DeVito *et al.* [8] that shows that a minimizer $f_{P,\lambda}$ of (3) exists.

**Proposition 8.** *Let $L$ be a convex loss function of $(a,p)$ type, let $P$ be a distribution on $X \times Y$ with $|P|_a < \infty$, let $H$ be a RKHS of a bounded kernel $k$ and let $\lambda > 0$. Then there exists a unique minimizer $f_{P,\lambda} \in H$ of $f \mapsto \mathcal{R}^{\mathrm{reg}}_{L,P,\lambda}(f)$ and $\|f_{P,\lambda}\|_H \leq \sqrt{\mathcal{R}_{L,P}(0)/\lambda} := \delta_{P,\lambda}$.*

If $H$ is a RKHS of a bounded kernel and $L$ is a convex, invariant loss of lower and upper order $p$, then it is easy to see by Lemma 7 that *exactly* for the distributions P with $|P|_p < \infty$, the minimizer $f_{P,\lambda}$ is uniquely determined. If $|P|_p = \infty$, we have $\mathcal{R}^{\mathrm{reg}}_{L,P,\lambda}(f) = \infty$ for *all* $f \in H$. Hence, we will use the definition $f_{P,\lambda} := 0$ for such P.

Our next aim is to establish a representation of $f_{P,\lambda}$. To this end, we define, for $p \in [1,\infty]$, the conjugate $p' \in [1,\infty]$ by $1/p + 1/p' = 1$. Furthermore, we have to recall the notion of subdifferentials; see Phelps [17].

**Definition 9.** *Let $H$ be a Hilbert space, let $F: H \to \mathbb{R} \cup \{\infty\}$ be a convex function and let $w \in H$ with $F(w) \neq \infty$. Then the subdifferential of $F$ at $w$ is defined by*

$$\partial F(w) := \{w^* \in H : \langle w^*, v - w \rangle \leq F(v) - F(w) \text{ for all } v \in H\}.$$

With the help of the subdifferential, we can now state the following theorem, which combines results from Zhang [30], Steinwart [22] and DeVito *et al.* [8].

**Theorem 10.** *Let $p \geq 1$, let $L$ be a convex loss function of type $(a,p)$ and let $P$ be a distribution on $X \times Y$ with $|P|_a < \infty$. Let $H$ be the RKHS of a bounded, continuous kernel $k$ over $X$ and let $\Phi: X \to H$ be the canonical feature map of $H$. Then there exists an $h \in L_{p'}(P)$ such that $h(x,y) \in \partial_2 L(y, f_{P,\lambda}(x))$ for all $(x,y) \in X \times Y$ and*

$$f_{P,\lambda} = -(2\lambda)^{-1} \mathbb{E}_P h \Phi, \tag{5}$$

*where $\partial_2 L$ denotes the subdifferential with respect to the second variable of $L$. Moreover, for all distributions $Q$ on $X \times Y$ with $|Q|_a < \infty$, we have $h \in L_{p'}(P) \cap L_1(Q)$ and*

$$\|f_{P,\lambda} - f_{Q,\lambda}\|_H \leq \lambda^{-1} \|\mathbb{E}_P h\Phi - \mathbb{E}_Q h\Phi\|_H, \tag{6}$$

*and if $L$ is an invariant loss of upper order $p$ and $|P|_p < \infty$, then $h \in L_{p'}(P) \cap L_{p'}(Q)$.*

## 3. Consistency of kernel-based regression

In this section we establish *L-risk consistency* of KBR methods, that is, we show that $\mathcal{R}_{L,P}(\hat{f}_{n,\lambda_n}) \to \mathcal{R}^*_{L,P}$ holds in probability for $n \to \infty$ and suitably chosen regularization sequences $(\lambda_n)$. Of course, such convergence can only hold if the RKHS is rich enough. One way to describe the richness of $H$ is the following definition taken from Steinwart [21].

***Definition 11.*** *Let $X \subset \mathbb{R}^d$ be compact and let $k$ be a continuous kernel on $X$. We say that $k$ is universal if its RKHS is dense in the space of continuous functions $(C(X), \|\cdot\|_\infty)$.*

It is well known that many popular kernels including the Gaussian RBF kernels are universal; see Steinwart [21] for a simple proof of the universality of the latter kernel. With Definition 11, we can now formulate our consistency result.

***Theorem 12.*** *Let $X \subset \mathbb{R}^d$ be compact, let $L$ be an invariant, convex loss of lower and upper order $p \geq 1$ and let $H$ be a RKHS of a universal kernel on $X$. Define $p^* := \max\{2p, p^2\}$ and fix a sequence $(\lambda_n) \subset (0, \infty)$ with $\lambda_n \to 0$ and $\lambda_n^{p^*} n \to \infty$. Then $\hat{f}_{n,\lambda_n}$ based on (2), using $\lambda_n$ for sample sets of length $n$, is $L$-risk consistent for all $P$ with $|P|_p < \infty$.*

Note that Theorem 12 in particular shows that KBR using the least squares loss function is *weakly universally consistent* in the sense of Györfi *et al.* [12]. Under the above assumptions on $L$, $H$ and $(\lambda_n)$, we can even *characterize* the distributions P for which KBR estimates based on (2) are $L$-risk consistent. Indeed, if $|P|_p = \infty$, then KBR is trivially $L$-risk consistent for P whenever $\mathcal{R}^*_{L,P} = \infty$. Conversely, if $|P|_p = \infty$ and $\mathcal{R}^*_{L,P} < \infty$, then KBR cannot be $L$-risk consistent for P because Lemma 7 shows $\mathcal{R}_{L,P}(f) = \infty$ for all $f \in H$.

In some sense it seems natural to consider only consistency for distributions that satisfy the tail assumption $|P|_p < \infty$, because this was done, for example, in Györfi *et al.* [12] for least squares methods. In this sense, Theorem 12 gives consistency for all reasonable distributions. However, the above characterization shows that our KBR methods are *not* robust against small violations of this tail assumption. Indeed, let P, and $\tilde{P}$ be two distributions with $|P|_p < \infty$, $|\tilde{P}|_p = \infty$ and $\mathcal{R}_{L,\tilde{P}}(f^*) < \infty$ for some $f^* \in L_p(P)$. Then every mixture distribution $Q_\varepsilon := (1-\varepsilon)P + \varepsilon\tilde{P}$, $\varepsilon \in (0,1)$, satisfies both $|Q_\varepsilon|_p = \infty$ and $\mathcal{R}^*_{L,Q_\varepsilon} < \infty$, and thus KBR is not consistent for any of the small perturbation $Q_\varepsilon$ of P, while it is consistent for original distribution P. From a robustness point of view, this is of course a negative result.

## 4. Robustness of kernel-based regression

In the statistical literature, different criteria have been proposed to define the notion of robustness in a mathematical way. In this paper, we mainly use the approach based on the influence function proposed by Hampel [13]. We consider a map $T$ that assigns to every distribution P on a given set $Z$, an element $T(P)$ of a given Banach space $E$. For the case of the convex risk minimization problem given in (3), we have $E = H$ and $T(P) = f_{P,\lambda}$. Denote the Dirac distribution at the point $z$ by $\Delta_z$, that is, $\Delta_z(\{z\}) = 1$.

***Definition 13*** (*Influence function*). *The influence function (IF) of $T$ at a point $z$ for a distribution P is the special Gâteaux derivative (if it exists)*

$$\mathrm{IF}(z; T, P) = \lim_{\varepsilon \downarrow 0} \varepsilon^{-1}\{T((1-\varepsilon)P + \varepsilon\Delta_z) - T(P)\}. \tag{7}$$

The influence function has the interpretation that it measures the impact of an infinitesimal small amount of contamination of the original distribution P in the direction of a Dirac distribution located in the point $z$ on the theoretical quantity of interest $T(\mathrm{P})$. Hence it is desirable that a statistical method $T(\mathrm{P})$ has a *bounded* influence function. We also use the sensitivity curve proposed by Tukey [27].

**Definition 14** (*Sensitivity curve*). *The sensitivity curve (SC) of an estimator $T_n$ at a point $z$ given a data set $z_1, \ldots, z_{n-1}$ is defined by*

$$\mathrm{SC}_n(z; T_n) = n(T_n(z_1, \ldots, z_{n-1}, z) - T_{n-1}(z_1, \ldots, z_{n-1})).$$

The sensitivity curve measures the impact of a single point $z$ and is a finite sample version of the influence function. If the estimator $T_n$ is defined via $T(\mathrm{P}_n)$, where $\mathrm{P}_n$ denotes the empirical distribution that corresponds to $z_1, \ldots, z_n$, then we have, for $\varepsilon_n = 1/n$,

$$\mathrm{SC}_n(z; T_n) = (T((1-\varepsilon_n)\mathrm{P}_{n-1} + \varepsilon_n \Delta_z) - T(\mathrm{P}_{n-1}))/\varepsilon_n. \tag{8}$$

In the following discussion, we give sufficient conditions for the existence of the influence function for the kernel-based regression methods based on (3). Furthermore, we establish conditions on the kernel $k$ and on the loss $L$ to ensure that the influence function and the sensitivity curve are bounded. Let us begin with the following results that ensure the existence of the influence function for KBR if the loss is convex and twice continuously differentiable.

**Theorem 15.** *Let $H$ be a RKHS of a bounded continuous kernel $k$ on $X$ with canonical feature map $\Phi: X \to H$ and let $L: Y \times \mathbb{R} \to [0, \infty)$ be a convex loss function of some strong type $(a, p)$. Furthermore, let P be a distribution on $X \times Y$ with $|\mathrm{P}|_a < \infty$. Then the influence function of $T(\mathrm{P}) := f_{\mathrm{P},\lambda}$ exists for all $z := (x, y) \in X \times Y$ and we have*

$$\mathrm{IF}(z; T, \mathrm{P}) = S^{-1}(\mathbb{E}_{\mathrm{P}}(L'(Y, f_{\mathrm{P},\lambda}(X))\Phi(X))) - L'(y, f_{\mathrm{P},\lambda}(x))S^{-1}\Phi(x), \tag{9}$$

*where $S: H \to H$, $S = 2\lambda \mathrm{id}_H + \mathbb{E}_{\mathrm{P}} L''(Y, f_{\mathrm{P},\lambda}(X))\langle \Phi(X), \cdot \rangle \Phi(X)$, is the Hessian of the regularized risk.*

It is worth mentioning that the proof can easily be modified to replace point mass contaminations $\Delta_z$ by arbitrary contaminations $\tilde{\mathrm{P}}$ that satisfy $|\tilde{\mathrm{P}}|_a < \infty$. As the discussion after Theorem 12 shows, we cannot omit this tail assumption on $\tilde{\mathrm{P}}$ in general.

From a robustness point of view, one is mainly interested in methods with bounded influence functions. Interestingly, for some kernel-based regression methods based on (3), Theorem 15 not only ensures the existence of the influence function, but also indicates how to guarantee its boundedness. Indeed, (9) shows that the only term of the influence function that depends on the point mass contamination $\Delta_z$ is

$$-L'(y, f_{\mathrm{P},\lambda}(x))S^{-1}\Phi(x). \tag{10}$$

Now recall that $\Phi$ is uniformly bounded because $k$ is assumed to be bounded. Consequently, the influence function is bounded in $z$ if and only if $-L'(y, f_{P,\lambda}(x))$ is uniformly bounded in $z = (x, y)$. Obviously, if, in addition, $L$ is invariant and $Y = \mathbb{R}$, then the latter condition is satisfied if and only if $L$ is Lipschitz continuous. Let us now assume that we use, for example, a Gaussian kernel on $X = \mathbb{R}^d$. Then the influence function is bounded *if and only if* $L$ is Lipschitz continuous. In particular, using the least squares loss in this scenario leads to a method with an unbounded influence function, while using the logistic loss function or its asymmetric generalization provides robust methods with bounded influence functions.

Unfortunately, the above results require a twice continuously differentiable loss and, therefore, they cannot be used, for example, to investigate methods based on the $\varepsilon$-insensitive loss or Huber's loss. Our next results, which in particular bound the difference quotient used in the definition of the influence function, apply to all convex loss functions of some type $(a, p)$ and hence partially resolve the above problem for non-differentiable losses.

**Theorem 16.** *Let $L : Y \times \mathbb{R} \to [0, \infty)$ be a convex loss of some type $(a, p)$, and let $P$ and $\tilde{P}$ be distributions on $X \times Y$ with $|P|_a < \infty$ and $|\tilde{P}|_a < \infty$. Furthermore, let $H$ be a RKHS of a bounded, continuous kernel on $X$. Then for all $\lambda > 0$, $\varepsilon > 0$, we have*

$$\|f_{(1-\varepsilon)P+\varepsilon\tilde{P},\lambda} - f_{P,\lambda}\|_H \leq 2c[\lambda\delta_{P,\lambda}]^{-1}\varepsilon(|P|_a + |\tilde{P}|_a + 2^{p+1}\delta_{P,\lambda}^p\|k\|_\infty^p + 2),$$

*where $c$ is the constant of the type $(a, p)$ inequality and $\delta_{P,\lambda} = \sqrt{\mathcal{R}_{L,P}(0)/\lambda}$.*

For the special case $\tilde{P} = \Delta_z$ with $z = (x, y)$, we have $|\Delta_z|_a = a(y)$ and hence we obtain bounds for the difference quotient that occurs in the definition of the influence function if we divide the bound by $\varepsilon$. Unfortunately, it then turns out that we can almost never bound the difference quotient *uniformly* in $z$ by the above result. The reason for this problem is that the $(a, p)$ type is a rather loose concept for describing the growth behaviour of loss functions. However, if we consider only *invariant* loss functions – and many loss functions used in practice are invariant – we are able to obtain stronger results.

**Theorem 17.** *Let $L : Y \times \mathbb{R} \to [0, \infty)$ be a convex invariant loss of upper order $p \geq 1$, and let $P$ and $\tilde{P}$ be distributions on $X \times Y$ with $|P|_p < \infty$ and $|\tilde{P}|_p < \infty$. Furthermore, let $H$ be a RKHS of a bounded, continuous kernel on $X$. Then for all $\lambda > 0$, $\varepsilon > 0$, we have*

$$\|f_{(1-\varepsilon)P+\varepsilon\tilde{P},\lambda} - f_{P,\lambda}\|_H \leq c\lambda^{-1}\|k\|_\infty\varepsilon|P - \tilde{P}|_{p-1} + |P - \tilde{P}|_0(\|k\|_\infty^{p-1}|P|_p^{(p-1)/2}\lambda^{(1-p)/2} + 1),$$

*where the constant $c$ only depends on $L$ and $p$. Moreover, if in addition $L$ is Lipschitz continuous, then for all $\lambda > 0$, $\varepsilon > 0$, we have*

$$\|f_{(1-\varepsilon)P+\varepsilon\tilde{P},\lambda} - f_{P,\lambda}\|_H \leq \lambda^{-1}\|k\|_\infty|L|_1\|P - \tilde{P}\|_\mathcal{M}\,\varepsilon.$$

*In particular considering (2), we have $\|\mathrm{SC}_n(z; T_n)\|_H \leq 2\lambda^{-1}\|k\|_\infty|L|_1$ for all $z \in X \times Y$.*

Finally, let us compare the influence function of kernel-based regression methods with the influence function of $M$ estimators in *linear* regression models with $f(x_i) = x_i'\theta$, where $\theta \in \mathbb{R}^d$ denotes the unknown parameter vector. Let us assume for reasons of simplicity that the scale parameter $\sigma \in (0, \infty)$ of the linear regression model is known. For more details about such $M$ estimators, see Hampel *et al.* [14]. The functional $T(\mathrm{P})$ that corresponds to an $M$ estimator is the solution of

$$\mathbb{E}_{\mathrm{P}} \eta(X, (Y - X'T(\mathrm{P}))/\sigma)X = 0, \tag{11}$$

where the odd function $\eta(x, \cdot)$ is continuous for $x \in \mathbb{R}^d$ and $\eta(x, u) \geq 0$ for all $x \in \mathbb{R}^d$, $u \in [0, \infty)$. Almost all proposals of $\eta$ may be written in the form $\eta(x, u) = \psi(v(x) \cdot u) \cdot w(x)$, where $\psi : \mathbb{R} \to \mathbb{R}$ is a suitable user-defined function (often continuous, bounded and increasing), and $w : \mathbb{R}^d \to [0, \infty)$ and $v : \mathbb{R}^d \to [0, \infty)$ are weight functions. An important subclass of $M$ estimators is of Mallows type, that is, $\eta(x, u) = \psi(u) \cdot w(x)$. The influence function of $T(\mathrm{P}) = \theta$ in the point $z = (x, y)$ at a distribution P for $(X, Y)$ on $\mathbb{R}^d \times \mathbb{R}$ is given by

$$\mathrm{IF}(z; T, \mathrm{P}) = M^{-1}(\eta, \mathrm{P}) \cdot \eta(x, (y - x'T(\mathrm{P}))/\sigma) \cdot x, \tag{12}$$

where $M(\eta, \mathrm{P}) := \mathbb{E}_{\mathrm{P}} \eta'(X, (Y - X'T(\mathrm{P}))/\sigma)XX'$. An important difference between kernel-based regression and $M$ estimation is that $\mathrm{IF}(z; T, \mathrm{P}) \in \mathbb{R}^d$ in (12), but $\mathrm{IF}(z; T, \mathrm{P}) \in H$ in (9) for point mass contamination in the point $z$.

A comparison of the influence function of KBR given in (9) with the influence function of $M$ estimators given in (12) yields that both influence functions have, nevertheless, a similar structure. The function $S = S(L'', \mathrm{P}, k)$ for KBR and the matrix $M(\eta, \mathrm{P})$ for $M$ estimation do not depend on $z$. The terms in the influence functions that depend on $z = (x, y)$, where the point mass contamination $\Delta_z$ occurs, are a product of two factors. The first factors are $-L'(y, f_{\mathrm{P}, \lambda}(x))$ for general KBR, $\psi(v(x) \cdot (y - x'\theta)/\sigma)$ for general $M$ estimation, $l'(y - f_{\mathrm{P}, \lambda}(x))$ for KBR with an invariant loss function and $\psi((y - x'\theta)/\sigma)$ for $M$ estimation of Mallows type. Hence the first factors measure the outlyingness in the $y$ direction. The KBR with an invariant loss function and $M$ estimators of Mallows type use first factors that depend only on the residuals. The second factors are $S^{-1}\Phi(x)$ for the kernel-based methods and $w(x)x$ for $M$ estimation. Therefore, they do not depend on $y$ and they measure the outlyingness in the $x$ direction.

In conclusion, one can say that there is a natural connection between KBR estimation and $M$ estimation in the sense of the influence function approach. The main difference between the influence functions is, of course, that the map $S^{-1}\Phi(x)$ takes values in the RKHS $H$ in the case of KBR, whereas $w(x)x \in \mathbb{R}^d$ for $M$ estimation.

## 5. Examples

In this section we give simple numerical examples to show the following concepts:

- The KBR with the $\varepsilon$-insensitive loss function is indeed more robust than the KBR based on the least squares loss function if there are outliers in the $y$ direction.

- In general, there is no hope to obtain robust predictions $\hat{f}(x)$ with KBR if $x$ belongs to a subset of the design space $X$ where no or almost no data points are in the training data set, that is, if $x$ is a leverage point.

We constructed a data set with $n = 101$ points in the following way. There is one explanatory variable $x_i$ with values from $-5$ to $5$ in steps of order $0.1$. The responses $y_i$ are simulated by $y_i = x_i + e_i$, where $e_i$ is a random number from a normal distribution with expectation 0 and variance 1. As hyperparameters, we used $(\varepsilon, \gamma, \lambda) = (0.1, 0.1, 0.05)$ for the $\varepsilon$-insensitive loss function with a RBF kernel, $(\varepsilon, \lambda) = (0.1, 0.05)$ for the $\varepsilon$-insensitive loss function with a linear kernel and $(\gamma, \lambda) = (0.1, 0.05)$ for the least squares loss function with an RBF kernel.

The $\varepsilon$-insensitive support vector regression ($\varepsilon$-SVR) and least squares support vector regression (LS-SVR) with similar hyperparameters give almost the same fitted curves; see Figure 1(a). However, Figure 1(b) shows that $\varepsilon$-SVR is much less influenced by outliers in the $y$ direction (one data point is moved to $(x, y) = (-2, 100)$) than is LS-SVR (cf. Wahba [29]) due to the different behaviour of the first derivative of the losses.

Now we add sequentially to the original data set three samples all equal to $(x, y) = (100, 0)$ which are bad leverage points with respect to a *linear* regression model. The number of such samples has a large impact on KBR with a linear kernel, but the predictions of KBR with a Gaussian RBF kernel are stable (but nonlinear), see Figure 1(c).

Now we study the impact of adding to the original data set two data points $z_1 = (100, 100)$ and $z_2 = (0, 100)$ on the predictions of KBR; see Figure 1(d). By construction, $z_1$ is a good leverage point and $z_2$ is a bad leverage point with respect to a *linear* regression model which follows, for example, by computing the highly robust least trimmed squares (LTS) estimator (Rousseeuw [19]), whereas the roles of these data points are switched for a quadratic model. There is no regression model which can fit all data points well because the $x$ components of $z_1$ and $z_2$ are equal by construction. We used the $\varepsilon$-insensitive loss function with a RBF kernel with hyperparameters $(\varepsilon, \gamma, \lambda) = (0.1, 0.1, 0.05)$ for the curve RBF, a and $(\varepsilon, \gamma, \lambda) = (0.1, 0.00001, 0.00005)$ for the curve RBF, b. This toy example shows that, in general, one cannot hope to obtain robust predictions $\hat{f}(x)$ for $\mathbb{E}_P(Y|X = x)$ with KBR if $x$ belongs to a subset of $X$ where no or almost no data points are in the training data set, because the addition of a single data point can have a big impact on KBR if the RKHS $H$ is rich enough. Note that the hyperparameters $\varepsilon$, $\gamma$ and $\lambda$ were specified in these examples to illustrate certain aspects of KBR and, hence, were not determined by a grid search or by cross-validation.

## 6. Discussion

In this paper, properties of kernel-based regression methods including support vector machines were investigated. Consistency of kernel-based regression methods was derived and results for the influence function, its difference quotient and the sensitivity curve were established. Our theoretical results show that KBR methods using a loss function with bounded first derivative (e.g., logistic loss) in combination with a bounded and rich enough continuous kernel (e.g., a Gaussian RBF kernel) are not only consistent and computational tractable, but also offer attractive robustness properties.

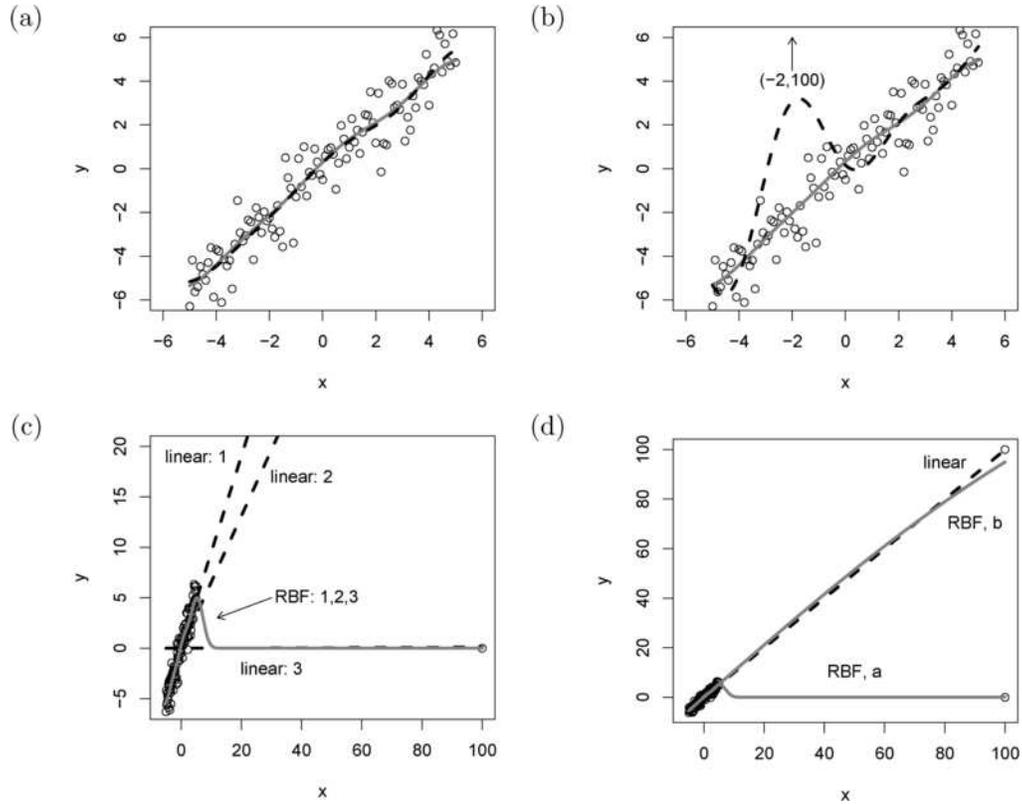

**Figure 1.** Results for simulated data sets without or with artificial outliers. (a) Linear relationship, no outliers: $\varepsilon$-insensitive loss function $L_\varepsilon$ (solid) and the least squares loss function $L_{LS}$ (dashed) both with a RBF kernel give almost the same results. (b) Linear relationship, one outlier in the $y$ direction at $(x,y) = (-2, 100)$: KBR with a RBF kernel performs more robustly if $L_\varepsilon$ (solid) is used instead of $L_{LS}$ (dashed). (c) Linear relationship with additional 1, 2 or 3 extreme points in $(x,y) = (100, 0)$: KBR using $L_\varepsilon$ with a linear kernel (dashed) and an RBF kernel (solid). (d) Linear relationship with two additional data points in $(x,y) = (100, 0)$ and $(x,y) = (100, 100)$: KBR with a linear kernel (dashed) and two curves based on RBF kernel (solid) with different values of $\lambda$ and $\gamma$.

Most of our results have analogues in the theory of kernel-based classification methods; see, for example, Christmann and Steinwart [5] and Steinwart [23]. However, because in the classification scenario $Y$ is only $\{-1, 1\}$ valued, many effects of the regression scenario with unbounded $Y$ do not occur in the above papers. Consequently, we had to develop a variety of new techniques and concepts: one central issue here was to find notions for loss functions which, on the one hand, are mild enough to cover a wide range of reasonable loss functions and, on the other hand, are strong enough to allow meaningful results for both consistency and robustness under minimal conditions on $Y$. In our analysis it turned out that the relationship between the growth behavior of the loss function and the tail behavior of $Y$ plays a central role for both types of results. Interestingly, similar tail properties of $Y$ are widely used to obtain consistency of nonparametric regression estimators and to establish robustness properties of $M$ estimators in linear regression. For example, Györfi *et al.* [12] assumed $\mathbb{E}_P Y^2 < \infty$ for the least squares loss, Hampel *et al.* ([14], page 315) assumed existence and non-singularity of $\mathbb{E}_P \eta'(X, (Y - X'T(P))/\sigma) X X'$ and Davies ([7], page 1876) assumed $\mathbb{E}_P \|X\|(\|X\| + |Y|) < \infty$. Another important issue was to deal with the estimation error in the consistency analysis. We decided to use a stability approach to avoid truncation techniques, so the proof of our consistency result became surprisingly short. An additional benefit of this approach was that it revealed an interesting connection between the robustness and the consistency of KBR methods. A somewhat similar observation was recently made by Poggio *et al.* [18] for a wide class of learning algorithms. However, they assumed that the loss or $Y$ is *bounded* and hence their results cannot be used in our more general setting.

Our result concerning the influence function of kernel-based regression (Theorem 15) are valid under the assumption that the loss function is twice continuously differentiable, whereas our other robustness results are valid for more general loss functions. The strong differentiability assumption was made because our proof is based on a classical theorem of implicit functions. We have not investigated whether similar results hold true for continuous but not differentiable loss functions. It may be possible to obtain such results by using an implicit function theorem for non-smooth functions based on a weaker concept than Fréchet differentiability. However, there are indications why a smooth loss function may even be desirable. The function $-L'$ has a role for kernel-based regression similar to the $\psi$ function for $M$ estimators. Huber ([16], page 51) considered robust estimation in parametric models and investigated the case that the underlying distribution is a mixture of a smooth distribution and a point mass. He showed that an $M$ estimator has a non-normal limiting behaviour if the point mass is at a discontinuity of the derivative of the score function. Because distributions with point masses are not excluded by nonparametric regression methods such as KBR, his results indicate that a twice continuously differentiable loss may guard against such phenomena.

To the best of our knowledge there are no results on robustness properties of KBR which are comparable to those presented here. However, we would like to refer to Schölkopf and Smola [20], who already gave arguments for better robustness properties of KBR when Huber's loss function is used instead of the least squares loss function.

Theorems 15 to 17 and the comments after Theorem 12 show that KBR estimators based on appropriate choices of $L$ and $k$ have a bounded influence function if the distribution P has the tail property $|P|_a < \infty$, but are non-robust against small violations of this tail assumption. The deeper reason for this instability is that the theoretical regularized risk itself is defined via $\mathbb{E}_P L(Y, f(X))$, which is a *non-robust* location estimator for the distribution of the losses. This location estimator can be infinite for mixture distributions $(1 - \varepsilon)P + \varepsilon \tilde{P}$ no matter how small $\varepsilon > 0$ is. Following general rules of robust estimation in linear regression models, one might replace this non-robust location estimator by a robust alternative like an $\alpha$-trimmed mean or the median (Rousseeuw [19]), which results in the kernel-based least median of squares estimator $f^*_{P,\lambda} = \arg\min_{f \in H} \text{Median}_P L(Y, f(X)) + \lambda\|f\|^2_H$. We conjecture that $f^*_{P,\lambda}$ offers additional robustness, but sacrifices computational efficiency. However, such methods are beyond the scope of this paper.

## Appendix: Proofs of the results

**Proof of Lemma 4.** The left inequality of (iii) is well known from convex analysis and the right inequality of (iii) easily follows from $l(r) = |l(r) - l(0)| \leq |r|V(r)$ for all $r \in \mathbb{R}$. Moreover, the right inequality of (iii) directly implies (ii). Furthermore, (i) can be easily obtained by the left inequality of (iii) because $\lim_{|r|\to\infty} l(r) = \infty$ implies $V(r) > 0$ for all $r \neq 0$ and the convexity of $l$ ensures that $r \mapsto l(|r|)$ is monotone. Finally, the last assertion follows from $L(y,t) = l(y-t) \leq c(|y-t|^p + 1) \leq c(|y|^p + |t|^p + 1)$. $\square$

**Proof of Proposition 6.** For bounded measurable functions $f: X \to \mathbb{R}$, we have $\mathcal{R}_{L,P}(f) \leq c\mathbb{E}_P(a(Y) + |f(X)|^p + 1) \leq c\|a\|_{L_1(P)} + c\|f\|^p_{L_p(P)} + c < \infty$. $\square$

**Proof of Lemma 7.** For all $a, b \in \mathbb{R}$, we have $(|a| + |b|)^p \leq 2^{p-1}(|a|^p + |b|^p)$ if $p \geq 1$ and $(|a| + |b|)^p \leq |a|^p + |b|^p$ otherwise. This obviously implies $|a|^p \leq 2^{p-1}(|a-b|^p + |b|^p)$ and $|a|^p \leq |a-b|^p + |b|^p$, respectively. If $f \in L_p(P)$, we consequently obtain

$$c\mathbb{E}_P(|Y|^p - c_p|f(X)|^p - 1) \leq c\mathbb{E}_P(|Y - f(X)|^p - 1) \leq \mathcal{R}_{L,P}(f) < \infty$$

for some finite constants $c > 0$ and $c_p > 0$. From this we immediately obtain $|P|_p < \infty$. The converse implication can be shown analogously. $\square$

**Proof of Proposition 8.** The first two assertions follow from Proposition 8 of DeVito *et al.* [8] and the last assertion is trivial. $\square$

**Proof of Theorem 10.** The existence of $h$ and (5) already have been shown by DeVito *et al.* [8]. Moreover, for $(x,y) \in X \times Y$ and $K_\lambda := \delta_{P,\lambda}\|k\|_\infty$, we have

$$|h(x,y)| \leq |L(y,.)|_{|[-f_{P,\lambda}(x), f_{P,\lambda}(x)]}|_1 \leq 2K_\lambda^{-1}\|L(y,.)|_{|[-4K_\lambda, 4K_\lambda]}\|_\infty$$
$$\leq 2cK_\lambda^{-1}(a(y) + |4K_\lambda|^p + 1) \qquad (13)$$

from which we deduce $h \in L_1(Q)$. Furthermore, by the definition of the subdifferential, we have $h(x,y)(f_{Q,\lambda}(x) - f_{P,\lambda}(x)) \leq L(y, f_{Q,\lambda}(x)) - L(y, f_{P,\lambda}(x))$ and hence

$$\mathbb{E}_{(X,Y)\sim Q}L(Y, f_{P,\lambda}(X)) + \langle f_{Q,\lambda} - f_{P,\lambda}, \mathbb{E}_Q h\Phi\rangle \leq \mathbb{E}_{(X,Y)\sim Q}L(Y, f_{Q,\lambda}(X)). \qquad (14)$$

Moreover, an easy calculation shows

$$\lambda\|f_{P,\lambda}\|^2_H + \langle f_{Q,\lambda} - f_{P,\lambda}, 2\lambda f_{P,\lambda}\rangle + \lambda\|f_{P,\lambda} - f_{Q,\lambda}\|^2_H = \lambda\|f_{Q,\lambda}\|^2_H. \qquad (15)$$

Combining (14), (15), $2\lambda f_{P,\lambda} = -\mathbb{E}_P h\Phi$ and $\mathcal{R}^{\text{reg}}_{L,Q,\lambda}(f_{Q,\lambda}) \leq \mathcal{R}^{\text{reg}}_{L,Q,\lambda}(f_{P,\lambda})$, it follows that

$$\mathcal{R}^{\text{reg}}_{L,Q,\lambda}(f_{P,\lambda}) + \langle f_{Q,\lambda} - f_{P,\lambda}, \mathbb{E}_Q h\Phi - \mathbb{E}_P h\Phi \rangle + \lambda \|f_{P,\lambda} - f_{Q,\lambda}\|_H^2 \leq \mathcal{R}^{\text{reg}}_{L,Q,\lambda}(f_{P,\lambda}).$$

Hence we obtain $\lambda \|f_{P,\lambda} - f_{Q,\lambda}\|_H^2 \leq \langle f_{P,\lambda} - f_{Q,\lambda}, \mathbb{E}_Q h\Phi - \mathbb{E}_P h\Phi \rangle \leq \|f_{P,\lambda} - f_{Q,\lambda}\|_H \cdot \|\mathbb{E}_Q h\Phi - \mathbb{E}_P h\Phi\|_H$, which shows the assertion in the general case.

Now let us assume that $L$ is invariant. As usual, we denote the function that represents $L$ by $l:\mathbb{R} \to [0,\infty)$. Then we easily check that $h$ satisfies $h(x,y) \in -\partial l(y - f_{P,\lambda}(x))$ for all $(x,y) \in X \times Y$. Now for $p = 1$ we see by (iv) of Lemma 4 that $l$ is Lipschitz continuous and hence the function $V$ is constant. Using $|h(x,y)| \leq V(y - f_{P,\lambda}(x))$ (cf. Phelps [17], Proposition 1.11) we thus find $h \in L_\infty(Q)$, which is the assertion for $p = 1$. Therefore, let us finally consider the case $p > 1$. Then for $(x,y) \in X \times Y$ with $r := |y - f_{P,\lambda}(x)| \geq 1$, we have $|h(x,y)| \leq |\partial l(y - f_{P,\lambda}(x))| \leq V(r) \leq \frac{2}{r}\|l_{|[2r,2r]}\|_\infty \leq cr^{p-1}$ for a suitable constant $c > 0$. Furthermore, for $(x,y) \in X \times Y$ with $|y - f_{P,\lambda}(x)| \leq 1$, we have $|h(x,y)| \leq |\partial l(y - f_{P,\lambda}(x))| \leq V(y - f_{P,\lambda}(x)) \leq V(1)$. Together, these estimates show $|h(x,y)| \leq \tilde{c}\max\{1, |y - f_{P,\lambda}(x)|^{p-1}\} \leq \tilde{c}\hat{c}_p(1 + |y|^{p-1} + |f_{P,\lambda}(x)|^{p-1})$ for some constant $\tilde{c}$ depending only on $L$ and $\hat{c}_p := \max\{1, 2^{p-2}\}$. Now, using $p'(p-1) = p$, we obtain

$$\|h\|_{L_{p'}(Q)} \leq \tilde{c}\hat{c}_p(|Q|_p + \|k\|_\infty^{p-1}\|f_{P,\lambda}\|_H^{p-1} + 1). \tag{16}$$

Finally, for later purpose we note that our previous considerations for $p = 1$ showed that (16) also holds in this case. □

To prove Theorem 12, we need the following preliminary results.

**Lemma 18.** *Suppose that the minimizer $f_{P,\lambda}$ of (3) exists for all $\lambda > 0$. Then we have*

$$\lim_{\lambda \to 0} \mathcal{R}^{\text{reg}}_{L,P,\lambda}(f_{P,\lambda}) = \inf_{f \in H} \mathcal{R}_{L,P}(f) := \mathcal{R}^*_{L,P,H}.$$

**Proof.** Let $\varepsilon > 0$ and $f_\varepsilon \in H$ with $\mathcal{R}_{L,P}(f_\varepsilon) \leq \mathcal{R}^*_{L,P,H} + \varepsilon$. Then for all $\lambda < \varepsilon \|f_\varepsilon\|_H^{-2}$, we have $\mathcal{R}^*_{L,P,H} \leq \lambda \|f_{P,\lambda}\|_H^2 + \mathcal{R}_{L,P}(f_{P,\lambda}) \leq \lambda \|f_\varepsilon\|_H^2 + \mathcal{R}_{L,P}(f_\varepsilon) \leq 2\varepsilon + \mathcal{R}^*_{L,P,H}$. □

**Lemma 19.** *Let $L$ be a convex and invariant loss function of lower and upper order $p \geq 1$ and let $H$ be a RKHS of a universal kernel. Then for all distributions $P$ on $X \times Y$ with $|P|_p < \infty$, we have $\mathcal{R}^*_{L,P,H} = \mathcal{R}^*_{L,P}$.*

**Proof.** Follows from Corollary 1 of Steinwart *et al.* [25]. □

**Lemma 20.** *Let $L$ be a convex invariant loss function of some type $p \geq 1$ and let $P$ be a distribution on $X \times Y$ with $|P|_p < \infty$. Then there exists a constant $c_p > 0$ that depends only on $L$ and $p$ such that for all bounded measurable functions $f, g: X \to Y$, we have*

$$|\mathcal{R}_{L,P}(f) - \mathcal{R}_{L,P}(g)| \leq c_p(|P|_{p-1} + \|f\|_\infty^{p-1} + \|g\|_\infty^{p-1} + 1)\|f - g\|_\infty.$$

**Proof.** Again we have $V(|y| + |a|) \leq \tilde{c}_p(|y|^{p-1} + |a|^{p-1} + 1)$ for all $a \in \mathbb{R}$, $y \in Y$, and a suitable constant $\tilde{c}_p > 0$ that depends on $L$ and $p$. Furthermore, we find

$$|\mathcal{R}_{L,P}(f) - \mathcal{R}_{L,P}(g)| \leq \mathbb{E}_P |l(Y - f(X)) - l(Y - g(X))|$$
$$\leq \mathbb{E}_P V(|Y| + \|f\|_\infty + \|g\|_\infty)|f(X) - g(X)|.$$

Now we easily obtain the assertion by combining both estimates. $\square$

**Lemma 21.** *Let $Z$ be a measurable space, let $P$ be a distribution on $Z$, let $H$ be a Hilbert space and let $g: Z \to H$ be a measurable function with $\|g\|_q := (\mathbb{E}_P \|g\|_H^q)^{1/q} < \infty$ for some $q \in (1, \infty)$. We write $q^* := \min\{1/2, 1/q'\}$. Then there exists a universal constant $c_q > 0$ such that, for all $\varepsilon > 0$ and all $n \geq 1$, we have*

$$P^n\left((z_1, \ldots, z_n) \in Z^n : \left\|\frac{1}{n}\sum_{i=1}^n g(z_i) - \mathbb{E}_P g\right\| \geq \varepsilon\right) \leq c_q \left(\frac{\|g\|_q}{\varepsilon n^{q^*}}\right)^q.$$

For the proof of Lemma 21, we have to recall some basics from local Banach space theory. To this end, we call a sequence of independent, symmetric $\{-1, +1\}$-valued random variables $(\varepsilon_i)$ a *Rademacher sequence*. Now let $E$ be a Banach space, let $(X_i)$ be an i.i.d. sequence of $E$-valued, centered random variables and let $(\varepsilon_i)$ be a Rademacher sequence which is independent of $(X_i)$. The distribution of $\varepsilon_i$ is denoted by $\nu$. Using Hoffmann-Jørgensen ([15], Corolllary 4.2), we have for all $1 \leq p < \infty$ and all $n \geq 1$ that

$$\mathbb{E}_{P^n}\left\|\sum_{i=1}^n X_i\right\|^p \leq 2^p \mathbb{E}_{P^n} \mathbb{E}_{\nu^n}\left\|\sum_{i=1}^n \varepsilon_i X_i\right\|^p, \tag{17}$$

where the left expectation is with respect to the distribution $P^n$ of $X_1, \ldots, X_n$, whereas the right expectation is also with respect to the distribution $\nu^n$ of $\varepsilon_1, \ldots, \varepsilon_n$. Furthermore, a Banach space $E$ is said to have type $p$, $1 \leq p \leq 2$, if there exists a constant $c_p(E) > 0$ such that for all $n \geq 1$ and all finite sequence $x_1, \ldots, x_n \in E$ we have $\mathbb{E}_{\nu^n}\|\sum_{i=1}^n \varepsilon_i x_i\|^p \leq c_p(E) \sum_{i=1}^n \|x_i\|^p$. In the following, because we are only interested in Hilbert spaces $H$, we note that these spaces always have type 2 with constant $c_2(H) = 1$ by orthogonality. Furthermore, they also have type $p$ for all $1 \leq p < 2$ by Kahane's inequality (see, e.g., Diestel *et al.* [10], page 211), which ensures $(\mathbb{E}_{\nu^n}\|\sum_{i=1}^n \varepsilon_i x_i\|^p)^{1/p} \leq c_{p,q}(\mathbb{E}_{\nu^n}\|\sum_{i=1}^n \varepsilon_i x_i\|^q)^{1/q}$ for all $p, q \in (0, \infty)$, $n \geq 1$, all Banach spaces $E$, all $x_1, \ldots, x_n \in E$ and constants $c_{p,q}$ only depending on $p$ and $q$.

**Proof of Lemma 21.** The summations used in this proof are with respect to $i \in \{1, \ldots, n\}$. Define $h(z_1, \ldots, z_n) := \frac{1}{n}\sum g(z_i) - \mathbb{E}_P g$. A standard calculation shows $P^n(\|h\| \geq \varepsilon) \leq \varepsilon^{-q} \mathbb{E}_{P^n}\|h\|_H^q$; hence it remains to estimate $\mathbb{E}_{P^n}\|h\|_H^q$. By (17) we have

$$\mathbb{E}_{P^n}\left\|\sum g(Z_i) - \mathbb{E}_P g\right\|_H^q \leq 2^q \mathbb{E}_{P^n} \mathbb{E}_{\nu^n}\left\|\sum \varepsilon_i (g(Z_i) - \mathbb{E}_P g)\right\|_H^q. \tag{18}$$

For $1 < q \leq 2$, we hence obtain

$$\mathbb{E}_{\mathrm{P}^n}\|h\|^q \leq 2^q n^{-q} \mathbb{E}_{\mathrm{P}^n}\mathbb{E}_{\nu^n}\left\|\sum \varepsilon_i(g(Z_i) - \mathbb{E}_\mathrm{P} g)\right\|^q \leq 2^q c_q n^{-q} \sum \mathbb{E}_\mathrm{P}\|g(Z_i) - \mathbb{E}_\mathrm{P} g\|^q$$
$$\leq 4^q c_q n^{1-q} \mathbb{E}_\mathrm{P}\|g\|^q,$$

where $c_q$ is the type $q$ constant of Hilbert spaces. From this we easily obtain the assertion for $1 < q \leq 2$. Now let us assume that $2 < q < \infty$. Then using (18) and Kahane's inequality there is a universal constant $c_q > 0$ with

$$\mathbb{E}_{\mathrm{P}^n}\|h\|_H^q \leq 2^q n^{-q} \mathbb{E}_{\mathrm{P}^n}\mathbb{E}_{\nu^n}\left\|\sum \varepsilon_i(g(Z_i) - \mathbb{E}_\mathrm{P} g)\right\|_H^q$$
$$\leq c_q n^{-q} \mathbb{E}_{\mathrm{P}^n}\left(\mathbb{E}_{\nu^n}\left\|\sum \varepsilon_i(g(Z_i) - \mathbb{E}_\mathrm{P} g)\right\|_H^2\right)^{q/2}$$
$$\leq c_q n^{-q} \mathbb{E}_{\mathrm{P}^n}\left(\sum \|g(Z_i) - \mathbb{E}_\mathrm{P} g\|_H^2\right)^{q/2}$$
$$\leq c_q n^{-q}\left(\sum (\mathbb{E}_\mathrm{P}\|g(Z_i) - \mathbb{E}_\mathrm{P} g\|_H^q)^{2/q}\right)^{q/2}$$
$$\leq 2^q c_q n^{-q/2} \mathbb{E}_\mathrm{P}\|g\|_H^q,$$

where in the third step we used that Hilbert spaces have type 2 with constant 1. From this estimate we easily obtain the assertion for $2 < q < \infty$. $\square$

**Proof of Theorem 12.** To avoid handling too many constants, we assume $\|k\|_\infty = 1$, $|\mathrm{P}|_p = 1$, and $c = 2^{-(p+2)}$ for the upper order constant of $L$. Then an easy calculation shows $\mathcal{R}_{L,\mathrm{P}}(0) \leq 1$. Furthermore, we assume without loss of generality that $\lambda_n \leq 1$ for all $n \geq 1$. This implies $\|f_{\mathrm{P},\lambda_n}\|_\infty \leq \|f_{\mathrm{P},\lambda_n}\|_H \leq \lambda_n^{-1/2}$. Now for $n \in \mathcal{N}$ and regularization parameter $\lambda_n$, let $h_n : X \times Y \to \mathbb{R}$ be the function obtained by Theorem 10. Then our assumptions and (16) give $\|h_n\|_{L_{p'}(\mathrm{P})} \leq 3 \cdot 2^{p^*/p-2} \lambda_n^{-(p-1)/2}$. Moreover, for $g \in H$ with $\|f_{\mathrm{P},\lambda_n} - g\|_H \leq 1$ we have $\|g\|_\infty \leq \|f_{\mathrm{P},\lambda_n}\|_\infty + \|f_{\mathrm{P},\lambda_n} - g\|_\infty \leq 2\lambda_n^{-1/2}$ and, hence, Lemma 20 provides a constant $c_p > 0$ that depends only on $p$ and $L$ such that

$$|\mathcal{R}_{L,\mathrm{P}}(f_{\mathrm{P},\lambda_n}) - \mathcal{R}_{L,\mathrm{P}}(g)| \leq c_p \lambda_n^{(1-p)/2}\|f_{\mathrm{P},\lambda_n} - g\|_H, \tag{19}$$

for all $g \in H$ with $\|f_{\mathrm{P},\lambda_n} - g\|_H \leq 1$. Now let $0 < \varepsilon \leq 1$ and $D \in (X \times Y)^n$ with corresponding empirical distribution D such that

$$\|\mathbb{E}_\mathrm{P} h_n \Phi - \mathbb{E}_\mathrm{D} h_n \Phi\|_H \leq c_p^{-1} \lambda_n^{(p+1)/2} \varepsilon. \tag{20}$$

Then Theorem 10 gives $\|f_{\mathrm{P},\lambda_n} - f_{\mathrm{D},\lambda_n}\|_H \leq c_p^{-1} \lambda_n^{(p-1)/2}\varepsilon \leq 1$ and, hence, (19) yields

$$|\mathcal{R}_{L,\mathrm{P}}(f_{\mathrm{P},\lambda_n}) - \mathcal{R}_{L,\mathrm{P}}(f_{\mathrm{D},\lambda_n})| \leq c_p \lambda_n^{-(p-1)/2}\|f_{\mathrm{P},\lambda_n} - f_{\mathrm{D},\lambda_n}\|_H \leq \varepsilon. \tag{21}$$

Let us now estimate the probability of $D$ satisfying (20). To this end we define $q := p'$ if $p > 1$ and $q := 2$ if $p = 1$. Then we have $q^* := \min\{1/2, 1/q'\} = \min\{1/2, 1/p\} = p/p^*$, and by Lemma 21 and $\|h_n\|_{L_{p'}(P)} \leq 3 \cdot 2^{p^*/p-2} \lambda_n^{-(p-1)/2}$ we obtain

$$P^n(D \in (X \times Y)^n : \|\mathbb{E}_P h_n \Phi - \mathbb{E}_D h_n \Phi\| \leq c_p^{-1} \lambda_n^{(p+1)/2} \varepsilon) \geq 1 - \hat{c}_p \left(\frac{3 \cdot 2^{p^*/p-2}}{\varepsilon \lambda_n^p n^{p/p^*}}\right)^{p'},$$

where $\hat{c}_p$ is a constant that depends only on $L$ and $p$. Now using $\lambda_n^p n^{p/p^*} = (\lambda_n^{p^*} n)^{p/p*} \to \infty$, we find that the probability of samples sets $D$ satisfying (20) converges to 1 if $n = |D| \to \infty$. As we have seen above, this implies that (21) holds true with probability tending to 1. Now, since $\lambda_n \to 0$ we additionally have $|\mathcal{R}_{L,P}(f_{P,\lambda_n}) - \mathcal{R}^*_{L,P}| \leq \varepsilon$ for all sufficiently large $n$ and hence we finally obtain the assertion. $\square$

To prove Theorem 15 we have to recall some notions from Banach space calculus. To this end, $B_E$ denotes the open unit ball of a Banach space $E$ throughout the rest of this section. We say that a map $G : E \to F$ between Banach spaces $E$ and $F$ is (Fréchet) differentiable in $x_0 \in E$ if there exists a bounded linear operator $A : E \to F$ and a function $\varphi : E \to F$ with $\frac{\varphi(x)}{\|x\|} \to 0$ for $x \to 0$ such that

$$G(x_0 + x) - G(x_0) = Ax + \varphi(x) \tag{22}$$

for all $x \in E$. Furthermore, because $A$ is uniquely determined by (22), we write $G'(x) := \frac{\partial G}{\partial E}(x) := A$. The map $G$ is called continuously differentiable if the map $x \mapsto G'(x)$ exists on $E$ and is continuous. Analogously we define continuous differentiability on open subsets of $E$. Moreover, we need the following two theorems which can be found, for example, in Akerkar [1] and Cheney [3], respectively.

**Theorem 22** (*Implicit function theorem*). *Let $E$ and $F$ be Banach spaces, and let $G : E \times F \to F$ be a continuously differentiable map. Suppose that we have $(x_0, y_0) \in E \times F$ such that $G(x_0, y_0) = 0$ and $\frac{\partial G}{\partial F}(x_0, y_0)$ is invertible. Then there exists a $\delta > 0$ and a continuously differentiable map $f : x_0 + \delta B_E \to y_0 + \delta B_F$ such that for all $x \in x_0 + \delta B_E$, $y \in y_0 + \delta B_F$ we have $G(x, y) = 0$ if and only if $y = f(x)$. Moreover, the derivative of $f$ is given by $f'(x) = -(\frac{\partial G}{\partial F}(x, f(x)))^{-1} \frac{\partial G}{\partial E}(x, f(x))$.*

**Theorem 23** (*Fredholm alternative*). *Let $E$ be a Banach space and let $S : E \to E$ be a compact operator. Then $\mathrm{id}_E + S$ is surjective if and only if it is injective.*

**Proof of Theorem 15.** The key ingredient of our analysis is the map $G : \mathbb{R} \times H \to H$ defined by $G(\varepsilon, f) := 2\lambda f + \mathbb{E}_{(1-\varepsilon)P + \varepsilon \Delta_z} L'(Y, f(X)) \Phi(X)$ for all $\varepsilon \in \mathbb{R}$, $f \in H$. Let us first check that its definition makes sense. To this end, recall that every $f \in H$ is a bounded function because we assumed that $H$ has a bounded kernel $k$. As in the proof of Proposition 6, we then find $\mathbb{E}_P |L'(Y, f(X))| < \infty$ for all $f \in H$. Because the boundedness of $k$ also ensures that $\Phi$ is a bounded map, we then see that the $H$-valued expectation used in the definition of $G$ is defined for all $\varepsilon \in \mathbb{R}$ and all $f \in H$. (Note that for $\varepsilon \notin [0, 1]$ the $H$-valued expectation is with respect to a signed measure; (cf. Dudley [11].) Now for $\varepsilon \in [0,1]$ we obtain (see Christmann and Steinwart [5] for a detailed derivation)

$$G(\varepsilon, f) = \frac{\partial R^{\text{reg}}_{L,(1-\varepsilon)\text{P}+\varepsilon\Delta_z,\lambda}}{\partial H}(f). \qquad (23)$$

Given an $\varepsilon \in [0,1]$, the map $f \mapsto R^{\text{reg}}_{L,(1-\varepsilon)\text{P}+\varepsilon\Delta_z,\lambda}(f)$ is convex and continuous by Lemma 20, and hence (23) shows that $G(\varepsilon, f) = 0$ if and only if $f = f_{(1-\varepsilon)\text{P}+\varepsilon\Delta_z,\lambda}$. Our aim is to show the existence of a differentiable function $\varepsilon \mapsto f_\varepsilon$ defined on a small interval $(-\delta, \delta)$ for some $\delta > 0$ that satisfies $G(\varepsilon, f_\varepsilon) = 0$ for all $\varepsilon \in (-\delta, \delta)$. Once we have shown the existence of this function, we immediately obtain $\text{IF}(z; T, \text{P}) = \frac{\partial f_\varepsilon}{\partial \varepsilon}(0)$. For the existence of $\varepsilon \mapsto f_\varepsilon$ we have to check by Theorem 22 that $G$ is continuously differentiable and that $\frac{\partial G}{\partial H}(0, f_{\text{P},\lambda})$ is invertible. Let us start with the first. By an easy calculation we obtain

$$\frac{\partial G}{\partial \varepsilon}(\varepsilon, f) = -\mathbb{E}_\text{P} L'(Y, f(X))\Phi(X) + \mathbb{E}_{\Delta_z} L'(Y, f(X))\Phi(X) \qquad (24)$$

and a slightly more involved computation (cf. Christmann and Steinwart [5]) shows

$$\frac{\partial G}{\partial H}(\varepsilon, f) = 2\lambda \, \text{id}_H + \mathbb{E}_{(1-\varepsilon)\text{P}+\varepsilon\Delta_z} L''(Y, f(X))\langle \Phi(X), \cdot \rangle \Phi(X) = S. \qquad (25)$$

To prove that $\frac{\partial G}{\partial \varepsilon}$ is continuous, we fix an $\varepsilon$ and a convergent sequence $f_n \to f$ in $H$. Because $H$ has a bounded kernel, the sequence of functions $(f_n)$ is then uniformly bounded. By the continuity of $L'$ we thus find a measurable bounded function $g : Y \to \mathbb{R}$ with $L'(y, f_n(x)) \leq L'(y, g(y))$ for all $n \geq 1$ and all $(x, y) \in X \times Y$. As in the proof of Proposition 6, we find $(y \mapsto L(y, g(y))) \in L_1(\text{P})$ and, therefore, an application of Lebesgue's theorem for Bochner integrals gives the continuity of $\frac{\partial G}{\partial \varepsilon}$. Because the continuity of $G$ and $\frac{\partial G}{\partial H}$ can be shown analogously, we obtain that $G$ is continuously differentiable (cf. Akerkar [1]). To show that $\frac{\partial G}{\partial H}(0, f_{\text{P},\lambda})$ is invertible it suffices to show by the Fredholm alternative that $\frac{\partial G}{\partial H}(0, f_{\text{P},\lambda})$ is injective and that $Ag := \mathbb{E}_\text{P} L''(Y, f_{\text{P},\lambda}(X))g(X)\Phi(X), g \in H$, defines a compact operator on $H$. To show the compactness of the operator $A$, recall that $X$ and $Y$ are Polish spaces (see Dudley [11]), because we assumed that $X$ and $Y$ are closed. Furthermore, Borel probability measures on Polish spaces are *regular* by Ulam's theorem, that is, they can be approximated from inside by compact sets. In our situation, this means that for all $n \geq 1$ there exists a compact subset $X_n \times Y_n \subset X \times Y$ with $\text{P}(X_n \times Y_n) \geq 1 - \frac{1}{n}$. Now we define a sequence of operators $A_n : H \to H$ by

$$A_n g := \int_{X_n} \int_{Y_n} L''(y, f_{\text{P},\lambda}(x)) \text{P}(\text{d}y|x) g(x) \Phi(x) \, \text{d}\text{P}_X(x) \qquad (26)$$

for all $g \in H$. Let us show that $A_n$, $n \geq 1$, is a compact operator. To this end we may assume without loss of generality that $\|k\|_\infty \leq 1$. For $g \in B_H$ and $x \in X$, we then have

$$h_g(x) := \int_{Y_n} L''(y, f_{\text{P},\lambda}(x))|g(x)|\text{P}(\text{d}y|x) \leq c \int_{Y_n} (a(y) + |f_{\text{P},\lambda}(x)|^p + 1)\text{P}(\text{d}y|x) =: h(x)$$

for $a:Y \to \mathbb{R}$, $p \geq 1$ and $c > 0$ according to the $(a,p)$ type of $L''$. Obviously, we have $\|h_g\|_1 \leq \|h\|_1 < \infty$ for all $g \in B_H$, and, consequently, $d\mu_g := h_g \, dP_X$ and $d\mu := h \, dP_X$ are finite measures. By Diestel and Uhl ([9], Corollary 8, page 48) we hence obtain

$$A_n g = \int_{X_n} \operatorname{sign} g(x) \Phi(x) \, d\mu_g(x) \in \mu(X_n) \overline{\operatorname{aco} \Phi(X_n)},$$

where $\operatorname{aco} \Phi(X_n)$ denotes the absolute convex hull of $\Phi(X_n)$, and the closure is with respect to $\|\cdot\|_H$. Now using the continuity of $\Phi$ we see that $\Phi(X_n)$ is compact and, hence, so is the closure of $\operatorname{aco} \Phi(X_n)$. This shows that $A_n$ is a compact operator. To see that $A$ is compact, it therefore suffices to show $\|A_n - A\| \to 0$ for $n \to \infty$. Define $B := (X \times Y) \setminus (X_n \times Y_n)$. Recalling that the convexity of $L$ implies $L'' \geq 0$, the latter convergence follows from $P(X_n \times Y_n) \geq 1 - \frac{1}{n}$, $L''(\cdot, f_{P,\lambda}(\cdot)) \in L_1(P)$ and

$$\|A_n g - A g\| \leq \int_B L''(y, f_{P,\lambda}(x)) |g(x)| \|\Phi(x)\| \, dP(x,y) \leq \|g\|_H \int_B L''(y, f_{P,\lambda}(x)) \, dP(x,y),$$

where again we assumed without loss of generality that $\|k\|_\infty \leq 1$. Let us now show that $\frac{\partial G}{\partial H}(0, f_{P,\lambda}) = 2\lambda \operatorname{id}_H + A$ is injective. To this end, let us choose $g \in H$ with $g \neq 0$. Then we find

$$\langle (2\lambda \operatorname{id}_H + A)g, (2\lambda \operatorname{id}_H + A)g \rangle > 4\lambda \langle g, Ag \rangle = 4\lambda \mathbb{E}_P L''(Y, f_{P,\lambda}(X)) g^2(X) \geq 0,$$

which shows the injectivity. As already described, we can now apply the implicit function theorem to see that $\varepsilon \mapsto f_\varepsilon$ is differentiable on $(-\delta, \delta)$. Furthermore, (24) and (25) yield

$$\operatorname{IF}(z; T, P) = \frac{\partial f_\varepsilon}{\partial \varepsilon}(0) = S^{-1}(\mathbb{E}_P(L'(Y, f_{P,\lambda}(X))\Phi(X))) - L'(y, f_{P,\lambda}(x)) S^{-1} \Phi(x).$$

$\square$

**Proof of Theorem 16.** By Theorem 10, there exists an $h \in L_\infty(P)$ such that

$$\|f_{P,\lambda} - f_{(1-\varepsilon)P + \varepsilon \tilde{P}, \lambda}\|_H \leq \varepsilon \lambda^{-1} \|\mathbb{E}_P h \Phi - \mathbb{E}_{\tilde{P}} h \Phi\|_H \leq \varepsilon \lambda^{-1} \|k\|_\infty (\|h\|_{L_1(P)} + \|h\|_{L_1(\tilde{P})}).$$

Because $h$ is independent of $\varepsilon$ and $\tilde{P}$, we then obtain the assertion by (13). $\square$

**Proof of Theorem 17.** Using the notation of the previous proof and $\|f_{P,\lambda}\|_\infty \leq \|k\|_\infty \times \sqrt{\mathcal{R}_{L,P}(0)/\lambda}$, we obtain analogously to the proof of Theorem 10 that

$$|h(x,y)| \leq \tilde{c}(|y|^{p-1} + |f_{P,\lambda}(x)|^{p-1} + 1) \leq c(|y|^{p-1} + \|k\|_\infty^{p-1} |P|_p^{(p-1)/2} \lambda^{-(p-1)/2} + 1),$$

where $\tilde{c}, c > 0$ are constants that depend only on $L$ and $p$. We then obtain the first assertion by combining the above estimate with $\|f_{P,\lambda} - f_{(1-\varepsilon)P + \varepsilon \tilde{P}, \lambda}\|_H \leq \varepsilon \lambda^{-1} \|\mathbb{E}_P h \Phi - \mathbb{E}_{\tilde{P}} h \Phi\|_H \leq \varepsilon \lambda^{-1} \|k\|_\infty \mathbb{E}_{|P - \tilde{P}|} |h|$. The second assertion can be shown analogously using $\|h\|_\infty \leq |L|_1$. Finally, the third assertion is a direct consequence of the second assertion and (8). $\square$

# Acknowledgements

We thank the editor and two anonymous referees for helpful comments.